\documentclass[12pt,reqno]{amsart}
\usepackage{amssymb,delarray}
\usepackage{amsfonts}
\usepackage{epsfig}
\usepackage[all]{xy}

\makeindex{}

\def\le{\leqslant}

\def\pr{{\sf pr}}
\def\pr{{\mbox{\rm pr}}}

\def\pr{\mathsf{pr}}

\newtheorem{thm}{Theorem}[section]
\newtheorem{lem}[thm]{Lemma}

\newtheorem{cor}[thm]{Corollary}

{\catcode`\@=11
\gdef\n@te#1#2{\leavevmode\vadjust{%
 {\setbox\z@\hbox to\z@{\strut#1}%
  \setbox\z@\hbox{\raise\dp\strutbox\box\z@}\ht\z@=\z@\dp\z@=\z@%
  #2\box\z@}}}
\gdef\leftnote#1{\n@te{\hss#1\quad}{}}
\gdef\rightnote#1{\n@te{\quad\kern-\leftskip#1\hss}{\moveright\hsize}}
\gdef\?{\FN@\qumark}
\gdef\qumark{\ifx\next"\DN@"##1"{\leftnote{\rm##1}}\else
 \DN@{\leftnote{\rm??}}\fi{\rm??}\next@}}

\begin{document}
\baselineskip=14.pt plus 2pt 

\title[Links and Hurwitz
curves]{Links and Hurwitz curves}
\author[Vik.S.~Kulikov]{Vik. S.~Kulikov}
\address{Steklov Mathematical Institute\\
Gubkina str., 8\\
119991 Moscow \\
Russia} \email{kulikov@mi.ras.ru}

\dedicatory{} \subjclass{}
\thanks{The work  was partially supported by
RFBR ({\rm No.} 02-01-00786). }
\begin{abstract}
In the note, we give a proof,  based on the Generalized Thom
Conjecture, of Bennequin's Theorem on upper bound for the Euler
number of a link which is considered as a closed braid. A lower
bound for the Euler number of a link is also given.
\end{abstract}

\maketitle
\setcounter{tocdepth}{2}


\def\st{{\sf st}}

\section{Introduction}

Let $l$ be a link in the three-dimensional sphere $S^3$ consisting
of $k$ components. Recall that an oriented  surface $S\subset S^3$
is called a {\it Seifert surface} of the link $l$ if the boundary
$\partial S$ of $S$ coincides with $l$ and $S$ has not a closed
component (without boundary). Let $\chi (S)$ be the Euler
characteristic of $S$. By definition, the {\it Euler number}
$e(l)$ of $l$ is
\begin{equation}
e(l)= \max _{S} \chi (S),
\end{equation}
where the maximum is taken over all Seifert surfaces of $l$. Note
that if $l$ is a knot of genus $g$, then
\begin{equation} e(l)=1-2g.
\end{equation}

By Alexander's theorem (see \cite{Al}), there is a number $m\in
\mathbb N$ such that a given link $l$ is equivalent to a closed
braid $\overline b$ (notation: $l\simeq \overline b$), where $b$
is a braid in the braid group $\mbox{\rm Br}_m$ on $m$ strings.

Below, we fix a set $\{ a_1,\dots ,a_{m-1} \} $ of so called {\it
standard generators} of $\mbox{\rm Br}_m$, i.e., generators being
subject to the relations
$$
\begin {array}{ll}
a_ia_{i+1}a_i & =a_{i+1}a_i a_{i+1} \qquad \qquad 1\leq i\leq m-2 ,  \\
a_ia_{k} & =a_{k}a_i  \qquad \qquad \qquad \, \, \mid i-k\mid \,
\geq 2
\end{array}
$$
and extend this set of generators to a set of generators $\{
a_{i,j}\}_{ 1\leq i< j\leq m } $, where $a_{i,i+1}= a_i$ and
\[
a_{i,j}= (a_{j-1}a_{j-2}\dots a_{i+1})a_i(a_{j-1}a_{j-2}\dots
a_{i+1})^{-1}
\]
for $j-i\geq 2$. An element $b\in \mbox{\rm Br}_m$ can be
presented as a word in the alphabet $\{ a_{i,j}, a_{i,j}^{-1} \}_{
1\leq i< j\leq m } $:
\begin{equation} \label{norm} b=w(a_{1,2},\dots,
a_{m-1,m})=\prod_{k=1}^{n_w} a_{i_k,j_k}^{\varepsilon_k},
\end{equation}
where $\varepsilon_k=\pm 1$. The minimum
\[ \mid b\mid = \min_{w(a_{i,j})=b} n_w, \]
where the minimum is taken over all presentations of $b$ in the
form (\ref{norm}) is called the {\it length} of $b$.

As is known, if braids $b_1$ and $b_2$ are conjugated in
$\mbox{\rm Br}_m$, then the closed braids $\overline b_1$ and
$\overline b_2$ are equivalent links. The number
\[ \mid\mid \overline b\mid \mid =\min_{g\in \mbox{\rm Br}_m}
 \mid g^{-1}bg \mid
\]
is called the {\it norm} of a closed braid $\overline b$.

Let $B_{l,m}= \{ b\in \mbox{\rm Br}_m\, \mid l\simeq \overline
b\}$ be the set of closed braids on $m$ strings equivalent to $l$.
If $B_{l,m}\neq \emptyset$, then the number
\[ \mid\mid l\mid \mid_m =\min_{b\in B_{l,m}}
\mid\mid \overline b \mid\mid
\]
is called the {\it $m$-norm} of a link $l$.

Denote by $\widetilde {\mbox{\rm Br}}^+_m$ the semigroup generated
 in the braid group
$\mbox{Br}_m$ by the set $\{ a_{i,j}\}_{ 1\leq i< j\leq m } $. An
element $b\in \mbox{\rm Br}_m$ is called {\it positive}
(respectively, {\it negative}) if $b\in \widetilde {\mbox{\rm
Br}}^+_m$ (respectively, if $b^{-1}\in \widetilde {\mbox{\rm
Br}}^+_m$).

Consider the homomorphism $\mbox{deg} : \mbox{\rm Br}_m \to
\mbox{Br}_m/[\mbox{\rm Br}_m,\mbox{\rm Br}_m ]\simeq \mathbb Z$
sending all $a_{i,j}$ to $1\in \mathbb Z$. The image $\mbox{\rm
deg}\, b$ of an element $b\in \mbox{\rm Br}_m$ is called the {\it
degree} of $b$.

The aim of this note is to give a proof, based on the Generalized
Thom Conjecture, of Bennequin's Theorem (\cite{Ben},\cite{Ben2})
on upper bound for the Euler number $e(l)$ in terms of invariants
of a closed braid $\overline b\simeq l$ and also to give some
lower bound for it.

\begin{thm}  \label{Be1} Let a link $l$ be presented as a closed braid
$\overline b$ for some $b\in \mbox{\rm Br}_m$. Then
\begin{equation} \label{main2}
m- \mid\mid \overline b \mid\mid \leq e(l).
\end{equation}
\end{thm}
\begin{thm} {\rm (\cite{Ben},\cite{Ben2})} \label{Be} Let a link $l$ be presented as a closed braid
$\overline b$ for some $b\in \mbox{\rm Br}_m$. Then
\begin{equation} \label{main}
e(l)\leq m-\mid \mbox{\rm deg}\, b\mid .
\end{equation}
\end{thm}

The idea of the proof of Theorem \ref{Be} is the following. First
of all, it is easy to see that the general case can be reduced to
the case $\deg \, b\geq 0$. Then for a given link $l\simeq
\overline b$, where $b\in \mbox{\rm Br}_m$, $\deg \, b\geq 0$,
applying results obtained in \cite{Kh-Ku} about so called Hurwitz
curves in the complex Hirzebruch surface $F_N$, we construct
smooth real surface $S$ and  algebraic curve $C$ lying in $F_N$
for some $N\geq 1$ and
having the genera $g(S)=1+
({Nm(m-1)-m -e(l)-\deg\, b})/{2}$ and $g(C)=1+(Nm(m-1)-2m)/2$, and
such that $[S]=[C]$, where $[C]$, $[S]$ $\in H_2(F_N,\mathbb Z)$
are the homology classes represented by real two-dimensional
surfaces $C$ and $S$. Now, the proof of Theorem \ref{Be} follows
from the Generalized Thom Conjecture proved in \cite{M-S-T} and
asserting that $g(C)\leq g(S)$.

 Since
$\mbox{\rm deg}\, b =\mid\mid \overline b \mid\mid $ for $b\in
\widetilde {\mbox{\rm Br}}^+_m$, we have the following corollary.
\begin{cor} Let a link $l$ be presented as a closed braid
$\overline b$ for some positive or  negative element $b\in
\mbox{\rm Br}_m$. Then
\begin{equation}
 \mid\mid l\mid \mid_m =\mid\mid \overline b \mid\mid=\mid
\mbox{\rm deg}\, b\mid ; \end{equation}
\begin{equation}
\label{main1}
 e(l)= m-\mid\mid l\mid \mid_m .
\end{equation}
\end{cor}

Obviously, $e(l)=k$ for a trivial link $l$ consisting of $k$
connected components. Therefore we have the following corollary.
\begin{cor} Let a link $l$ consisting of $k$ connected components
be presented as a closed braid
$\overline b$ for some element $b\in \mbox{\rm Br}_m$. If $$k>
m-\mid \mbox{\rm deg}\, b\mid $$ then $l$ is a non-trivial link.
\end{cor}

{\it Acknowledgement.} The author thanks I.A.Dynnikov for
references and his helpful remarks  during the preparation of this
paper.

\section{Proof of theorem \ref{Be1} }
To prove Theorem \ref{Be1}, let us identify the sphere $S^3$ with
the boundary $\partial D=(\partial D_1)\times D_2\cup D_1\times
\partial D_2$ of a bi-disc
\[
D=D_1\times D_2=\{ (z,w)\in \mathbb C^2\, \mid \, \, \mid
z\mid\leq 1,\, \, \mid w\mid \leq 2\}.
\]
Choose $m$ points $w_k=e^{\frac{2\pi  \sqrt{-1} k}{m}}\in D_2=\{
\mid w\mid \leq 2\}$, $k=1,\dots,m$, and identify the braid group
$\mbox{\rm Br}_m$ with the braid group $\mbox{\rm Br}[D_2,\{
w_1,\dots,w_m\} ]$. In this case the generators $a_{i,j}$ are
identified with half-twists along the segments $w=tw_i+(1-t)w_j$,
$t\in [0,1]$ (see Fig. 1), and $\overline b$ with a closed braid
lying in
$(\partial D_1)\times D_2$. 

\smallskip
\noindent\includegraphics[width=\hsize]{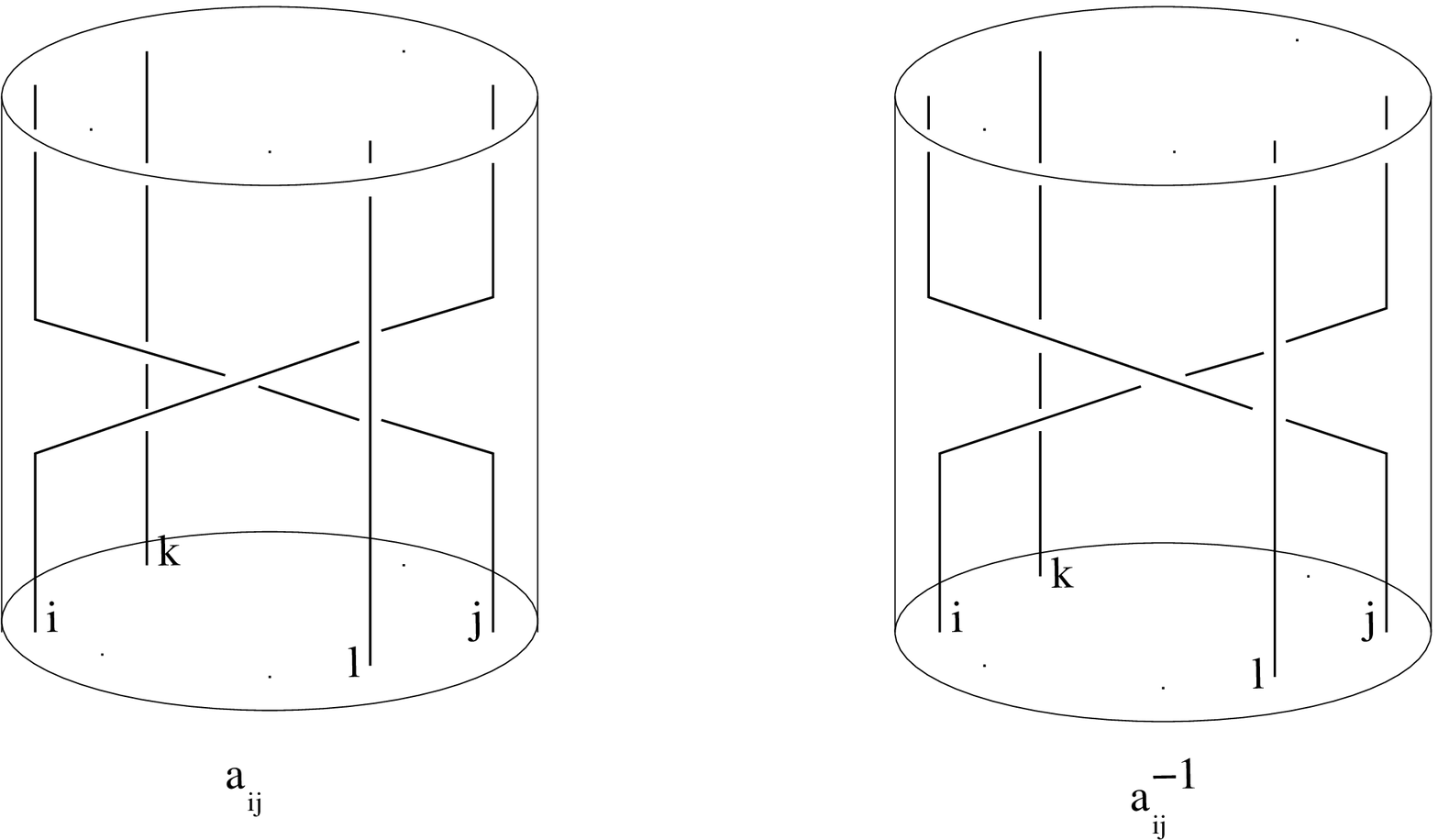}

\smallskip

\begin{center}
Fig. 1 \end{center}
 Let a link $l\simeq \overline b$, where
$$b=\prod_{k=1}^{n_b} a_{i_k,j_k}^{\varepsilon_k}\in \mbox{\rm
Br}_m,\, \, \, \varepsilon _k=\pm 1.$$ In this case one can
construct a Seifert surface $S$ of the link $l$ similar to the
construction in the standard case when the link $l$ is represented
as a projection of $l$ to a plane whose image is an immersed curve
with simple intersections (Wirtinger presentation). Namely, take
$m$ discs
$$S_j=\{(z,w)\in S^3\, \mid \,
 \mid z\mid\leq 1,\, \, w=2e^{\frac{2\pi  \sqrt{-1} j}{m}}\}\subset D_1\times \partial D_2, $$
 $j=1,\dots ,m$, glue each $S_j$ along a circle
 $$C_j=\{(z,w)\in S^3\, \mid \,
 \mid z\mid = 1,\, \, w=2e^{\frac{2\pi  \sqrt{-1} j}{m}}\} \subset \partial D_1\times \partial D_2$$
 with an annulus
$$A_j=\{(z,w)\in S^3\, \mid \,
 \mid z\mid = 1,\, \, w=2te^{\frac{2\pi  \sqrt{-1} j}{m}}+(1-t)e^{\frac{2\pi  \sqrt{-1} j}{m}},\, \, t\in [0,1] \} $$
 and put $\overline S_j= S_j\cup_{C_j}A_j$. Obviously, each $\overline
 S_j$ is a disc. Next, in each $$(\partial
D_1)_k\times D_2=\{ (z,w)\in (\partial D_1)\times D_2\, \mid \, \,
z=e^{\frac{2\pi  \sqrt{-1} t}{n_b}},\, \, k-\frac{1}{3}\leq t\leq
k+\frac{1}{3}  \}$$ let us attach a band  $B_k\simeq [0,1]\times
[0,1]$ to $\overline S_{i_k}$ and $\overline S_{j_k}$ in
dependence  on the sign of $\varepsilon _k$ as it is depicted in
Fig. 2.

 As a result, we obtain a surface $S$ in the sphere $S^3$
with the boundary $\overline b$. Obviously, the Euler
characteristic $\chi (S)=m-n_b$. Therefore,  Theorem \ref{Be1} is
proven.

\smallskip
\noindent\includegraphics[width=\hsize]{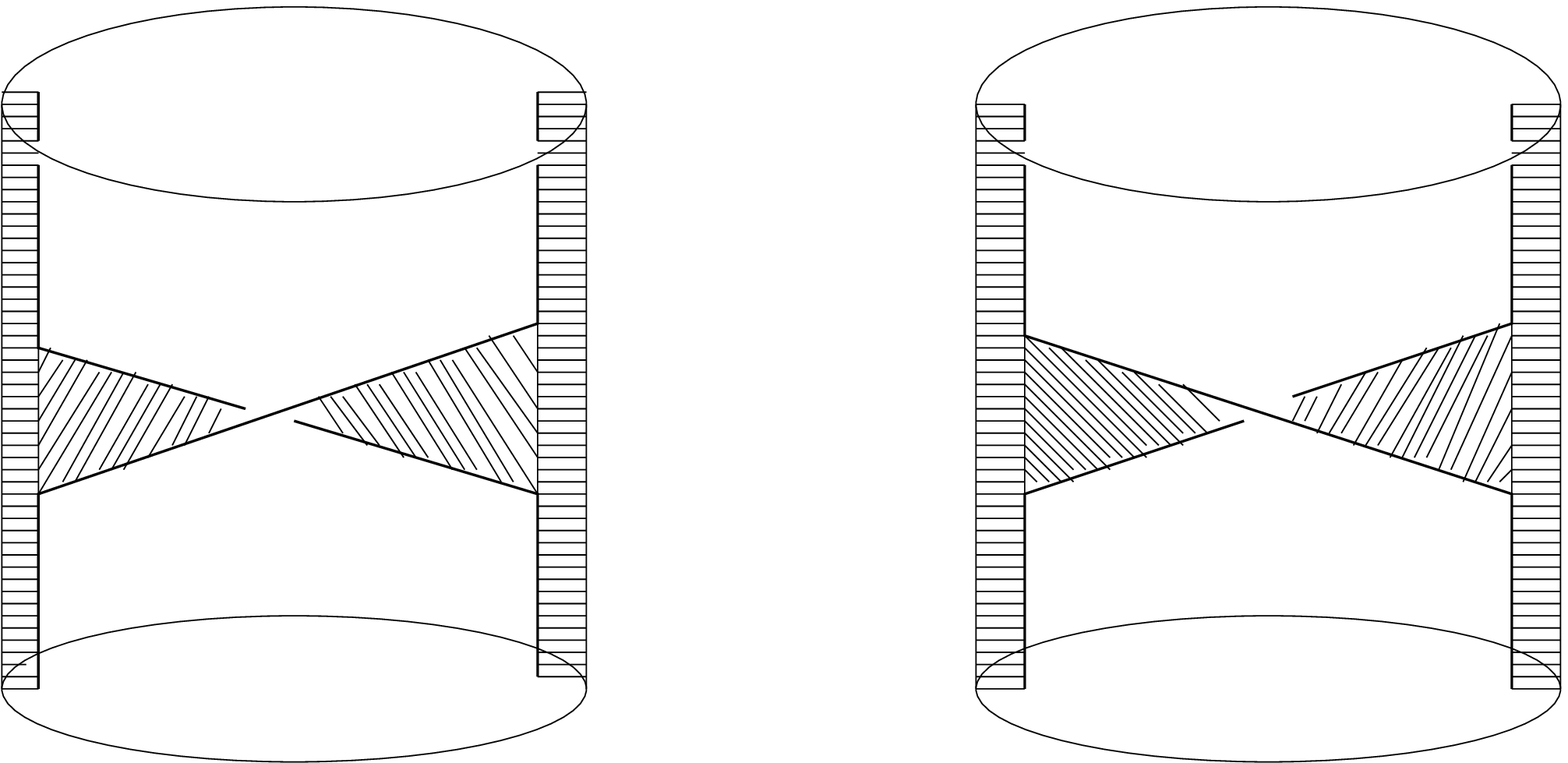}

\smallskip

\begin{center}
Fig. 2 \end{center}

\section{Proof of theorem \ref{Be} }
To prove Theorem \ref{Be}, let us, in the beginning, briefly
recall definitions of topological Hurwitz curves and their braid
monodromy factorizations given in \cite{Kh-Ku}. For a group
$\mbox{\rm Br}_m$ one can define the {\it factorization semigroup}
$S_{\mbox{\rm Br}_m}$. For this, consider an alphabet
$$X= \{ x_g\,  \mid\,  g\in \mbox{\rm Br}_m\}$$ and two sets of relations:
\newline $ R_{g_1,g_2;r}$ stands for  $x_{g_1}\cdot
x_{g_2}=x_{g_2}\cdot\,x_{g_2^{-1}g_1g_2}$ if $g_2\ne\bf 1$ and
$x_{g_1}\cdot\,x_{\bf 1}=x_{g_1}$; $R_{g_1,g_2;l}$ stands for
$x_{g_1}\cdot x_{g_2}=x_{g_1g_2g_1^{-1}}\cdot\,x_{g_1}$ if
$g_1\ne\bf 1$ and $x_{\bf 1}\cdot x_{g_2}=x_{g_2}$. Now, put
$$
\mathcal{R}=\{ R_{g_1,g_2;r},
 R_{g_1,g_2;l} \,\vert \, (x_{g_1},x_{g_2})\in X\times X,\, g_1\ne g_2\,
 \}
$$
and introduce the semigroup
$$S_{\mbox{\rm Br}_m}= \langle\, x\in X\,\,
:\,\, R\in \mathcal{R}\, \rangle 
$$
by means of this relation set $\mathcal R$. Introduce also a {\it
product homomorphism} $\alpha :S_{\mbox{\rm Br}_m}\to \mbox{\rm
Br}_m$ given by $\alpha (x_g)=g$ for each $x_{g}\in X$.

Denote by $F_N$  a relatively minimal ruled rational complex
surface (a {\it Hirzebruch surface}), $N\geq 1$, $\pr: F_N \to
\mathbb C\mathbb P^1$ the ruling, $R$ a fiber of $\pr$ and $E_N$
the exceptional section, $E_N^2=-N$. The variety $F_N\setminus
(E_N\cup R)$ is naturally isomorphic to  the complex affine plane
$\mathbb C^2$ with complex coordinates $(z,w)$ such that $\pr
(z,w)=z$.

By definition, the image $\bar H=f(S)\subset F_N$ of a continuous
map $f:S \to F_N\setminus E_N$ of an oriented closed real surface
$S$ is called a  {\it topological Hurwitz curve} (in $F_N$) of
{\it degree} $m$ if there is a finite subset $Z\subset\bar H$ such
that:
\begin{itemize}
 \item[(i)] $f$ is a smooth embedding
of the surface $S\setminus f^{-1}(Z)$ and for any $p\notin Z$,
$\bar H$ and the fiber $R_{\pr(p)}$ of $\pr$ meet at $p$
transversely with positive intersection number;
\item[(ii)] the
restriction of $\pr$ to $\bar H$ is a finite map of degree $m$.
{\rm (We call a map finite if the preimage of each point is
finite.) }
\end{itemize}
Choose a fibre $R=R_{\infty}$ being in general position with a
topological Hurwitz curve $\bar H$.  Put $\mathbb C^2=F_N\setminus
(E_N\cup R_{\infty})$ and fix complex coordinates $(z,w)$ in
$\mathbb C^2$ such that $\pr (z,w)=z$.
 At any point $p\in Z$ there is a well-defined ($W$-prepared) germ
 $(D,H=\bar H\cap D,\pr)$ of this curve in a bi-disc $D= D_1\times D_2$,
$D_1=D_1(\epsilon_1)=\{\mid z-z(p)\mid \leq \epsilon_1\}$,
$D_2=D_2(\epsilon_2)= \{ \mid w-w(p)\mid \leq \epsilon_2\}$,
$0<\epsilon_1<<\epsilon_2$, centered at $p$ and such that the
restriction of $\pr$ to $H$ is a proper map of a finite degree
$k\leq m$. If $\epsilon_1,\epsilon_2$ are sufficiently small,
then: $R_{\pr(p)}\cap H=p$; the above degree does not depend on
$\epsilon_1,\epsilon_2$; and the link $\partial D\cap H$ defines a
unique, up to conjugation, braid $b\in \mbox{\rm Br}_k\subset
\mbox{\rm Br}_m$, where $k$ is the above degree. So that, we may
speak on a {\it $tH$-singularity} $(D,H,\pr)$  of {\it degree} $k$
and {\it type} $b$.

When we are given a link $l\subset\partial D_1\times D_2$
realizing a braid $b\in B_k$, we associate with it a {\it standard
conical model} of a topological singularity of type $b$. It is
given by $H=C(l)$,
$$C(l)=\{ (rz,rw)\,\, \vert \,0\le r\le 1, (z,w)\in l \}. $$

As is known (see, for example, \cite{Kh-Ku}), if ($D,C,\pr$) is a
germ of a $W$-prepared $tH$-singularity then the germ ($D,C,\pr$)
is homeomorphic to the cone singularity of type
$b=\pr^{-1}(\partial D_1)\cap C$.

Since $\bar H\cap E_N=\emptyset$, one can define a {\it braid
monodromy factorization} $b(\bar H)\in S_{\mbox{\rm Br}_m}$ of
$\bar H$. For doing this, we fix a fiber $R_{\infty}$ meeting
transversely $\bar H$ and consider $\bar H\cap \mathbb C^2$, where
$\mathbb C^2=F_N\setminus (E_N\cup R_{\infty})$. Choose $r_1>>1$
such that $\pr (Z)\subset D_{1}(r_1)=\{ \mid z\mid \leq
r_1\}\subset\mathbb
C= \mathbb C \mathbb P^1\setminus\pr (R_\infty)$. 
Denote by $z_1,\dots,z_n$ the elements of the set $\pr (Z)$ and
assume that for each $i$ the intersection $\pr^{-1}(z_i)\cap Z$
consists of a single point. Pick $\rho$, $0<\rho<<1$, such that
the discs $D_{1,i}(\rho)=\{ z\in \mathbb C \, \mid \,\, \mid
z-z_i\mid < \rho \, \}$, $i=1,\dots, n$, would be disjoint. Select
arbitrary points $u_i\in
\partial D_{1,i}(\rho)$ and a point $u_0\in \partial D_{1}(r)$. Let
$D_{2}(r_2)=\{ w\in \mathbb C\, \, \mid\, \,  \, \mid w\mid \leq
r_2\}$ be a disc of radius $r_2>>1$ such that $\bar H\cap
\pr^{-1}(D_1(r_1))\subset D_1(r_1)\times D_2(r_2)$. Put
$D_{2,u_0}=\{ (u_0,w)\in \mathbb C^2\, \, \mid\, \,  \, \mid w\mid
\leq r_2\} \subset \pr^{-1}(u_0)$,
$K(u_0)=\{w_{1},\dots,w_{m}\}=D_{2,u_0}\cap \bar H$, and
$\mbox{\rm Br}_m=\mbox{\rm Br}[D_{2,u_0},K(u_0)]$. Choose disjoint
simple paths $l_i\subset \overline D_{1}(r_1)\setminus \bigcup^n_1
D_{1,i}(\rho)$, $i=1,\dots,n$, starting at $u_0$ and ending at
$u_i$ and renumber the points in a way that the product
$\gamma_1\dots\gamma_n$ of the loops $\gamma_i=l_i\circ \partial
D_{1,i}(\rho)\circ l_i^{-1}$ would be equal to $\partial D_1(r_1)$
in $\pi_1(\overline D_{1}(r_1)\setminus \{ z_1, \dots ,z_n\},
u_0)$. Each $\gamma _i$ defines an element $b_i\in {\mbox{\rm
Br}_m}$ represented by the paths $\pr^{-1}(\gamma_i)\cap \bar H$
starting and ending at the points lying in $K(u_0)$. The
factorization $b(\bar H)=x_{b_1}\cdot .\, .\, .\cdot x_{b_n}\in
S_{\mbox{\rm Br}_m}$ is called a {\it braid monodromy
factorization} of $\bar H$.

Denote by
$$\Delta^2_m=(a_{1,2}a_{2,3}\dots a_{m-1,m})^m$$
a generator of the center of  the group $\mbox{\rm Br}_m$. It is
easy to prove the following lemma (see, for example, \cite{Kh-Ku})
\begin{lem} For a topological Hurwitz curve $\bar H\subset F_N$ it holds
$$\alpha (b(\bar H))=\Delta^{2N}.$$
\end{lem}

The converse statement can be also proved in a straightforward
way.

\begin{thm} {\rm (}\cite{Kh-Ku}{\rm )} For any $s=x_{b_1}\cdot .\, .\, .\cdot
x_{b_n}\in S_{\mbox{\rm Br}_m}$ such that $\alpha (s)=\Delta^{2N}$
there is a topological Hurwitz curve $\bar H\subset F_N$ with a
braid monodromy factorization $b(\bar H)$ equal to $s$.
\end{thm}

Now we are able to prove  inequality  (\ref{main}). First of all,
it easy to see that if $l\simeq \overline b$ for some $b\in
\mbox{\rm Br}_m$, then the link $\overline {b^{-1}}$ is equivalent
to the mirror-image $\widetilde {l^{-1}}$ of the inverted link
$l^{-1}$. Therefore, to prove inequality (\ref{main}), we can
assume that $\deg\, b\geq 0$, since $e(l)=e(l^{-1})=e(\widetilde
{l^{-1}})$.

It follows from Theorem 5 in \cite{G} (see, for example, Lemma 1.3
in \cite{Kh-Ku}) that for any $b\in B_m$ there is a positive
element $r\in \widetilde {\mbox{\rm Br}}^+_m$ and a positive
integer $N\geq 1$ such that $rb=\Delta_m^{2N}$. We have $\deg
\Delta_m^2=m(m-1)$. Therefore $\deg\, r= Nm(m-1)-\deg\, b>0$. Let
\begin{equation} \label{r} r =\prod_{k=1}^{\deg\, r}
a_{i_k,j_k}
\end{equation}
be a presentation of $r$ as a word in the alphabet $\{
a_{i,j}\}_{1\leq i<j\leq m}$. Factorization (\ref{r}) defines an
element
$$s=( \prod_{k=1}^{\deg\, r}x_{a_{i_k,j_k}})\cdot x_b$$
in the factorization semigroup $S_{\mbox{\rm Br}_m}$. The element
$s$ is a braid monodromy factorization of a topological Hurwitz
curve $\overline H\subset F_N$ whose set $\pr (Z)$ of the critical
values  consists of points $z_k=\deg\, r-k+2$ for $k=1,\dots,
\deg\, r$ and $z_{\deg\, r+1}=0$, and whose braid monodromy  over
$z_k$, $k=1,\dots, \deg\, r$, is equal to the $k$-th factor
$a_{i_k,j_k}$ entering  in (\ref{r}), and whose braid monodromy
over the point $z_0$ is equal to $b$. Moreover, without loss of
generality, we can assume that $\overline H\cap \mbox{\rm
pr}^{-1}(\partial D_1)= \overline b\subset (\partial D_1)\times
D_2$, where $D_1=\{ \mid z\mid \leq 1\}$ and $D_2=\{ \mid w\mid
\leq r\}$ for some $r>>1$. Since all $a_{i,j}$ are conjugated to
$a_{1,2}$ and the element $a_{1,2}$ is the monodromy of the
critical value of the function given by $w^2=z$, then we can
assume that the Hurwitz curve $S_2=\overline H\cap \mbox{\rm
pr}^{-1}(\mathbb C\mathbb P^1\setminus D_1)$ is a smooth real
surface in $F_N$.

Consider the restriction of $\mbox{\rm pr}$ to $S_2$:
$$\mbox{\rm pr}_{\mid S_2}:S_2\to D_{\geq 1}=\mathbb C\mathbb P^1\setminus D_{1}.$$
The Euler characteristic of $S_2$ is equal to
$$\chi
(S_2)=m -Nm(m-1)+\deg\, b,$$
 since $D_{\geq 1}$ is a disc,
$\mbox{\rm pr}_{\mid S_2}$ has $\deg\, r= Nm(m-1)-\deg\, b$
simplest critical values, and $\deg \pr_{\mid S_2}=m$.

Let $S_1\subset \partial (D_1\times D_2)$ be a Seifert surface of
the link $\overline b\simeq l$. We can assume that $$\chi
(S_1)=e(l).$$ Consider a surface $S$ in $F_N$ which is obtained
from $S_1$ and $S_2$ by gluing along $\overline b$. Obviously, $S$
is a  closed real surface. Without loss of generality (after small
deformation of $S$ near $\overline b$), we can assume that $S$ is
a smooth surface. Since $\chi (\overline b)=0$, the Euler
characteristic
\begin{equation} \label{e1} \chi (S)=
\chi (S_1)+\chi (S_2)= e(l)+m -Nm(m-1)+\deg\, b.
\end{equation}

Consider the class $[S]$ of $S$ in the homology group
$H_2(F_N,\mathbb Z)$. As is known, the group $H_2(F_N,\mathbb Z)$
is generated by the class $[R]$ of a fibre $R$ of $\mbox{\rm pr}$
and the class $[E_N]$ of the exceptional section $E_N$ which have
the following intersection numbers: $[R]\cdot [R]=0$, $[R]\cdot
[E_N]=1$, and $[E_N]\cdot [E_N]=-N$. By construction of $S$, we
have $[S]\cdot [R]=\deg\, \bar H=m$ (to see this, one can consider
the intersection of $\overline H$ and a fibre $R_z$ lying over a
point $z\in D_{\geq 1}$) and $[S]\cdot [E_N]=0$, since $S_1\subset
D\subset \mathbb C^2\subset F_N\setminus E_N$ and by definition of
topological Hurwitz curves, $\overline H\cap E_N=\emptyset$.
Therefore
$$[S]=m[E_N]+Nm[R].$$

Let $C\subset F_N$ be a non-singular algebraic curve whose class
$[C]=mE_N+Nm[R]$. It is well-known that its genus
\begin{equation} \label{g2}
g(C)=(Nm(m-1)-2m)/2 +1.
\end{equation}

Since $C\subset F_N$ is an algebraic non-singular curve, it
follows from the Generalized Thom Conjecture proved in
\cite{M-S-T} that $\chi (S)\leq \chi (C)=2-2g(C)$ for any smooth
surface $S\subset F_N$ whose class $[S]=[C]$ in $H_2(F_N,\mathbb
Z)$. Therefore, applying (\ref{e1}) and (\ref{g2}), we have
$$\chi (S)=e(l)+m -Nm(m-1)+\deg\, b\leq 2m- Nm(m-1).$$
Thus,
$$e(l)\leq m-\deg\, b.$$

\ifx\undefined\bysame
\newcommand{\bysame}{\leavevmode\hbox to3em{\hrulefill}\,}
\fi

\end{document}